\documentclass[10pt]{amsart}
\usepackage[]{amsmath, amsthm, amsfonts}
\usepackage{graphicx}
\usepackage[all]{xy}

\newcommand {\C} {{\mathbb C}}
\newcommand {\R} {{\mathbb R}}
\newcommand {\Z} {{\mathbb Z}}
\newcommand {\Q} {{\mathbb Q}}

\newcommand {\dt} {\bullet}
\newcommand {\F} {{\mathcal F}}
\newcommand {\I} {{\mathcal I}}

\newcommand {\cH} {{\mathcal H}}
\newcommand {\G} {{\mathcal G}}
\newcommand {\LL} {{\mathcal L}}
\newcommand {\A} {{\mathcal A}}
\newcommand {\E} {{\mathcal E}}
\newcommand {\B} {{\mathcal B}}
\newcommand {\M} {{\mathcal M}}
\newcommand {\dd} {{\partial}}
\renewcommand {\AA} {{\mathbb A}}
\newcommand {\OO} {{\mathcal O}}

\newcommand {\Xan} {{X^{\text{an}}}}

\newtheorem{thm}[subsection]{Theorem}
\newtheorem{cor}[subsection]{Corollary}
\newtheorem{lemma}[subsection]{Lemma}
\newtheorem{prop}[subsection]{Proposition}

\newtheorem{ex}[subsection]{Example}

\begin{document}

\title{Mixed Hodge structures associated to geometric variations}

\author{ Donu Arapura}
 \address{
Department of Mathematics\\
 Purdue University\\
 West Lafayette, IN 47907\\
 U.S.A.}
 \thanks{Author partially supported by the NSF}

\email{dvb@math.purdue.edu}

\maketitle

In \cite{arapura}, I gave a construction
of a  mixed Hodge structure on the cohomology of a geometric
variation of Hodge structure as a subquotient of the   mixed Hodge
structure on the total space. I will refer to this as the
naive mixed Hodge structure.
One of the aims of this paper is to show that Morihiko Saito's theory of mixed
Hodge modules  yields the same mixed Hodge 
structure.  Although the proof of this equivalence can be
made quite short by simply quoting the basic properties of mixed Hodge
modules, I am not comfortable with doing this. 
Saito's theory has been laid out in a
pair of long and densely written papers \cite{saito1,saito2}, and it
would be fair to say that the details have not been very widely assimilated.
 Therefore I have decided to outline mixed Hodge module theory, to the extent that I need to,
 in order to describe Saito's mixed Hodge structure in more explicit
 terms. So to a large extent this paper is expository.

 A variation of Hodge structure gives a good notion of a regular family
 of Hodge structures. It is natural to extend this to families with singularities.
 Recall that  the constituents of a variation of a Hodge structure are a local
 system and a compatible filtered vector bundle with connection. As a
 first step toward the more general notion of a Hodge module, 
 these are replaced by a perverse sheaf and a filtered
 $D$-module respectively.  At this stage the challenge for any
 good theory would be to find a subclass of these
 pairs which is  both robust and meaningful  Hodge theoretically. Saito does this by a
 rather delicate induction. He   demands that the
 ``restriction'' of such a pair to any subvariety, given by a  vanishing cycle
 construction, should stay in the subclass; and that over a point it
determines a Hodge structure in the usual sense. For mixed Hodge
 modules, Saito also imposes  a weight filtration with additional
 constraints. The category of mixed Hodge modules possesses direct
 images. So in particular the cohomology of these modules carry mixed
 Hodge structures; among these are the Hodge structures indicated in the first paragraph.

 Here is a quick synopsis. The first three sections
summarize background material on the Riemann-Hilbert correspondence, 
vanishing cycles, and standard Hodge theory. All of these ingredients
are needed for mixed Hodge modules  discussed  in the fourth
section. The explicit description of Saito's mixed Hodge structure is given in 
\S\ref{section:expconstr}.
The comparison with the naive Hodge structure is contained in the
fifth section.
I have also included an appendix which clarifies the  proof of a lemma 
of \cite{arapura} which, as  Mark de Cataldo 
pointed out me, contained a hidden assumption. I would like to thank
to V. Srinivas for arranging a visit to  TIFR, and for suggesting 
that I give talks  on some of this material.

\section{D-modules and perverse sheaves}

\subsection{D-modules}
Suppose that  $X$ is a smooth complex algebraic variety.
 A $\C$-linear endomorphism of $\OO_X(U)$ such that
$[\ldots [[P,f_1],f_2],\ldots f_k]=0$ for all $f_i\in\OO_X(U)$ is
called an algebraic differential operator of order at most $k$ on
$U\subseteq X$. These form a sheaf of noncommutative rings $D_X$.
 When $X =\mathbb{A}^n$, the
ring of global sections of $D_X$ is the so called Weyl algebra, which is given
by generators
$x_1,\ldots x_n$, $\dd_1 \ldots \dd_n$ 
and relations
$$[x_i,x_j]= [\dd_i,\dd_j]=0$$
$$[\dd_i,x_j] = \delta_{ij}$$ 
Using local analytic coordinates, a similar description is available
in general. In other words, the stalk of the sheaf of holomorphic
differential operators $D_\Xan$ (defined as above) is isomorphic to the stalk of
$D_{\AA^{n.an}}$ at $0$. $D_X$ has an increasing filtration $F_kD_X$ consisting
of operators of order at most $k$. A standard computation
shows that the associated graded for the Weyl algebra
$$Gr (\Gamma(D_{\mathbb{A}^n})) = \bigoplus_k F_k/F_{k-1}$$
is just the polynomial ring in generators $x_1,\ldots x_n$, $\dd_1
\ldots \dd_n$. In general, $Gr(D_X)$ can be identified with the sheaf
of regular functions on the cotangent bundle $T^*X$ pushed down to $X$.

Since $D_X$ is noncommutative, we have to take care to
distinguish between left and right modules over it. We usually mean
left. We will primarily be interested in modules which are
coherent, i.e. locally finitely presented, over $D_X$.
An important example of a coherent $D_X$-module is given by a vector bundle
$V$ with an integrable algebraic connection $\nabla$. Pretty much by definition
the connection gives an action of first order differential operators
on $V$, and the  integrability condition
$$[\nabla_{\dd_i},\nabla_{\dd_j}] = 0$$
ensures that this extends to all of $D_X$. There are plenty of other
examples, such as $D_X$ itself, which do not arise this way.
A large class of examples can be given as follows.
 Let $X= \AA^n$ with coordinates $x_1\ldots x_n$, $Y = \AA^{n+1}$ with
 an additional coordinate
  $ x_{n+1}$ and suppose $i(x_1,\ldots, x_n) = (x_1,\ldots x_n,0)$. Given a
  $D_X$-module $M$, set
$$i_+M = \bigoplus_{j=0}^\infty\dd_{n+1}^j M$$
(using Borel's notation \cite{borel}) where  $\dd_{n+1}^j$ are treated as symbols. 
This becomes a $D_Y$-module via
$$
P\dd_{n+1}^jm\mapsto
\begin{cases}
  \dd_{n+1}^{j}Pm & \text{if $P\in \{x_1,\ldots ,x_n,  \dd_1, \ldots, \dd_n\}$}\\
0 & \text{if $P= x_{n+1}$}\\
\dd_{n+1}^{j+1}m & \text{if $P =\dd_{n+1}$}
\end{cases}
$$
The direct image $i_+M $ of a $D$-module along a more general closed immersion can be
defined by a similar procedure.

We define a {\em good } filtration on a $D_X$-module $M$
to be a filtration $F_p M$ such that
\begin{enumerate}
\item The filtration $F_pM= 0$ for $p\ll 0$ and $\cup F_pM = M$.
\item Each $F_pM$ is a coherent $\OO_X$-submodule.
\item $F_pD_X\cdot F_qM\subseteq F_{p+q}M$.
\end{enumerate}

\begin{thm}
Let $M$ be a nonzero coherent $D_X$-module. Then
  \begin{enumerate}
  \item   $M$ possess a good filtration.
\item The support of $Gr(M)$ in $T^*X$, called the characteristic variety, depends only on $M$.
\item (Bernstein's inequality) The dimension of the characteristic
  variety is greater than or equal to $\dim X$.
  \end{enumerate}
\end{thm}

A module $M$ is called {\em holonomic} if either it is zero or equality holds in (3).
The characteristic variety of an integrable connection is the zero
section of $T^*X$, so it is holonomic. Whereas $D_X$ is not holonomic
since its characteristic variety is all of $T^*X$. If $M$ is
holonomic, then so is $i_+M$ for any inclusion. This is easy to see if $i$ is
the inclusion of a point $p$, since the characteristic variety of
$i_+M$ is $T_p$. The definition of holonomicity 
given  is geometric, there is also a homological characterization
which yields:

\begin{thm}
  The full subcategory of holonomic modules is abelian and artinian
  (i.e. the descending chain condition holds).
\end{thm}

Recall that an algebraic integrable connection $\nabla$ on a variety $X$ has  regular
singularities if for suitable compactification, the connection matrix
has simple poles along the boundary divisor.
One can show that  a simple holonomic module always restricts to
an integrable connection on a nonempty open subset of its support.
This can be used to extend the above notion  
to holonomic modules. Namely a holonomic module is {\em regular} if
all of it simple subquotients restrict to connections with regular singularities.

Let $\Omega_\Xan^p$ denote the sheaf of holomorphic $p$-forms on the
associated complex manifold $\Xan$.  
We can modify the standard de Rham complex to allow coefficients in any
$D_\Xan$-module $M$:
$$DR(M)^\dt =  \Omega_\Xan^\dt\otimes_{\OO_\Xan} M[\dim X]$$
with differential  given by
$$d(dx_{i_1}\wedge\ldots \wedge d{x_{i_p}}\otimes m) = \sum_j dx_j\wedge
dx_{i_1}\wedge\ldots d{x_{i_p}}\otimes \dd_jm$$
The shift by $\dim X$ in $DR$ is done for convenience, as it simplifies various formulas. For
example, we have
$$DR(M)^\dt \cong \R {\mathcal H}om(\OO_{\Xan}, M)$$
which gives an extension   of $DR$ to the derived category $D^b(D_X)$.
When $M$ comes from an integrable connection $\nabla$, $DR(M)$ gives
a resolution of the locally constant sheaf $\ker(\nabla)$ (up to 
shift). Conversely any locally constant sheaf arises from such
a $D$-module. This is the classical version of the Riemann-Hilbert correspondence.
 For a general $M$, $DR(M)$ will no longer be a 
locally constant sheaf in general, but rather a complex with constructible
cohomology. Recall that  a  $\C_\Xan$-module $L$ is
{\em constructible} if there exists an algebraic stratification of $X$
such that the  restrictions $L$ to the strata are locally constant with finite dimensional
stalks.

\begin{thm}[Kashiwara, Mebkhout]
  The de Rham functor $DR$ induces an equivalence of categories
  between the subcategory $D^b_{rh}(D_X)\subset D^b(D_X)$ of complexes
  with regular holonomic cohomology
 and the  subcategory
$D_{constr}^b(\C_\Xan)\subset D^b(\C_\Xan)$
of complexes with constructible cohomology.
\end{thm}

There is more to the story. Various standard sheaf theoretical operations,
such as inverse and direct images,  correspond to natural
operations in the $D$-module world. See \cite{borel, kashiwara} for further details. 

\subsection{Perverse sheaves}
The category of  of regular holonomic modules sits in the triangulated category
$D_{rh}^b(D_X)$ as an abelian subcategory. Its image
 under $DR$ is the abelian category of  complex {\em perverse sheaves}
 $Perv(\C_{\Xan})$ \cite{bbd}.
In spite of the name, these objects are neither perverse nor sheaves,
but rather a class of well behaved elements of
$D_{constr}^b(\C_\Xan)$.

\begin{ex}
  $DR(\OO_X) = \C_X[\dim X]$ is perverse. More generally $L[\dim X]$
  is perverse for any local system.
\end{ex}

\begin{ex}
 Suppose that $X$ is complete (e.g. projective).
 Suppose that  $V$ is a vector bundle with   an integrable connection
$\nabla:V\to \Omega_X^1(\log D)\otimes V$ with logarithmic singularities along
a normal crossing divisor.
Let $U = X-D$ and $j:U\to X$ be the
inclusion. Then we have 
$$DR(j_*j^*V) = \R j_*(L)[\dim X]$$
is perverse. 
\end{ex}

\begin{ex}\label{ex:RjD}
  We have the $i$th perverse cohomology functor
  ${}^p\mathcal{H}^i:D^b_{constr}(X)\to Perv(\C)$ which corresponds to the
functor which assigns  the $i$th cohomology to a complex of
regular holonomic $D$-modules.
\end{ex}

Perverse sheaves can be characterized by purely sheaf theoretic
methods:
\begin{thm}
  $F\in D^b_{constr}(\C_\Xan)$ is perverse if and only if
\begin{enumerate}
\item For all $i$, $\dim supp\, \mathcal{H}^{-i}(F)\le i$.
\item These inequalities also hold for the Verdier dual 
$$D(F) = \R \mathcal{H}om(F, \C_{\Xan}[2\dim X])$$
\end{enumerate}
\end{thm}

These conditions can be verified directly in example \ref{ex:RjD}. 
The first condition follows from the fact that
$R^ij_*L$ is supported on the union of $i$-fold intersections of
components of $D$. For the second, it is enough to observe that
 $D (\R j_*L[\dim X]) = j_!L^\vee [\dim X]$ satisfies (1).
Note that the above conditions work perfectly well with other
coefficients, such as $\Q$, to define full subcategories $Perv(\Q_\Xan)\subset
D_{constr}^b(\Q_\Xan)$. Perverse sheaves have another source, independent of
$D$-modules.
In the 1970's Goresky and Macpherson
introduced intersection homology by a geometric construction
by placing restrictions how chains met the singular set in
terms of a function referred to as the  perversity.
 Their motivation was to find a theory which behaved like ordinary
homology for nonsingular spaces in general; for example, by satisfying
Poincar\'e duality.  When their constructions were
recast in sheaf theoretic
language, they provided 
basic examples of perverse sheaves.

\begin{ex}\label{ex:IC}
Suppose that  $Z\subseteq X$ is a possible singular subvariety.
Then the complex $IC_Z(\Q)$ computing the rational intersection
cohomology of $Z$ is (after a suitable shift and extension to $X$) 
a  perverse sheaf on $X$.
This is more generally true for the complex $IC_Z(L)$ computing
intersection cohomology of $Z$ with coefficients in
a locally constant sheaf defined on a Zariski
open $U\subseteq Z$.
   In the notation of
\cite{bbd}, this would be denoted by $i_*j_{!*}L[\dim Z]$, where
$j:U\to Z$ and $i:Z\to X$ are the inclusions. In the simplest
case, where $U\subseteq Z=X$ is the complement of a finite set,
$IC_X(\Q)= \tau_{<\dim X} \R j_*\Q[\dim X]$, where $\tau$ is the
truncation operator.
\end{ex}

The objects  $IC_Z(L)$ have a central place in the theory, since
all the simple objects are known to be of this form.
It follows that all perverse sheaves can
built up from such sheaves, since the category is artinian.
Further details can be found in \cite{bbd, brylinksi,goresky}.

\section{Vanishing cycles}

\subsection{Vanishing cycles}

Vanishing cycle sheaves and their corresponding $D$-modules
form the basis for Saito's constructions described later.
We will start with the classical picture.
 Suppose that $f:X\to \C$ is a morphism from a
nonsingular variety. 
The fiber $X_0=f^{-1}(0)$ may be singular, but the nearby 
fibers $X_t$, $0<|t|< \epsilon\ll 1$ are not. The premiage of
the $\epsilon$-disk $f^{-1}\Delta_\epsilon$ retracts onto $X_0$,
and $f^{-1}(\Delta_\epsilon-\{0\})\to \Delta_\epsilon-\{0\}$ is a 
fiber bundle. Thus we have a monodromy action by the (counterclockwise) generator $T\in
\pi_1(\C^*,t)$ on $H^i(X_t)$. (From now on, we will tend to treat algebraic varieties as an analytic
spaces, and will no longer be scrupulous about making a distinction.)
 The image of the restriction map
$$H^i(X_0)=H^i(f^{-1}\Delta_\epsilon)\to H^i(X_t),$$
lies in the kernel of $T-1$. The restriction is dual to the map
in homology which is induced by the (nonholomorphic)
 collapsing map of  $X_t$ onto $X_0$; the cycles which die in
the process are the vanishing cycles.

Let us reformulate things in a more abstract way following \cite{sga7}.
 The {\em nearby cycle} functor  applied to  $F\in D^b(X)$ is
$$\R \Psi F = i^{*} \R p_* p^* F,$$
where $\tilde \C^* $ is the universal cover of $\C^*=\C-\{0\}$, and $p:\tilde
\C^*\times_\C X\to X$, $i:X_0=f^{-1}(0)\to X$ are the natural maps. 
The {\em vanishing cycle} functor  $\R\Phi F$ is the mapping cone of the adjunction
morphism  $i^{*}F\to \R\Psi F$, and hence it fits into a distinguished triangle
$$i^*F\to \R\Psi F \stackrel{can}{\longrightarrow} \R\Phi F\to i^*F[1]$$
Both $\R \Psi F$ and $\R \Phi F$ are
often both loosely referred to as sheaves of vanishing cycles.
These objects possess  natural monodromy actions by $T$.
If we give $i^*F$ the trivial $T$ action, then
the  diagram with solid arrows commutes. 
$$
\xymatrix{
 i^*F\ar[r]\ar[d] & \R\psi_* F\ar[r]^{can}\ar[d]^{T-1} & \R\phi_*F\ar[r]\ar@{-->}[d]^{var} & i^*F[1]\ar[d] \\ 
 0\ar[r] & \R\psi_* F\ar[r]^{=} & \R\psi_* F\ar[r] & 0
}
$$
Thus we can deduce a morphism $var$, which completes this to a morphism
of triangles. In particular,  $T-1= var\circ can$.
 A further diagram chase also shows that $can \circ var = T-1$.

Given $p\in X_0$, let $B_\epsilon$ be an $\epsilon$-ball in $X$ centered
at $p$. Then $f^{-1}(t)\cap B_\epsilon$ is the so called Milnor fiber.
The stalks
$$\mathcal{H}^i(\R \Psi \Q)_p = H^i(f^{-1}(t)\cap B_\epsilon,\Q)$$
$$\mathcal{H}^i(\R \Phi \Q)_p = \tilde H^i(f^{-1}(t)\cap B_\epsilon,\Q)$$
give the (reduced) cohomology of the Milnor fiber. And
$$H^i(X_0, \R \Psi \Q) = H^i(f^{-1}(t),\Q)$$
is, as the terminology suggests, the cohomology of the nearby fiber.
We have a long exact sequence
$$\ldots H^i(X_0,\Q)\to  H^i(X_t,\Q)\stackrel{can}{\longrightarrow}
 H^i(X_0, \R \Phi\Q)\to\ldots$$

The following is a key ingredient in the whole story \cite{bbd,brylinksi}:
\begin{thm}[Gabber]
  If $L$ is perverse, then so are $\R\Psi L[-1]$ and $\R\Phi L [-1]$.
\end{thm}

We set ${}^p\psi_f L = {}^p\psi L = \R\Psi L[-1]$ and  
${}^p\phi_f L = {}^p\phi L = \R\Phi L[-1]$.

\subsection{Perverse Sheaves on a polydisk}\label{section:poldisk}

Let  $\Delta$ be a disk with the standard coordinate function $t$,
and inclusion  $j:\Delta-\{0\}=\Delta^*\to \Delta$. For simplicity assume $1\in \Delta^*$.
Consider a perverse sheaf $F$ on $\Delta$ which is locally constant
on $\Delta^*$. Then we can form the diagram
$$
\xymatrix{
 {}^p\psi_t F\ar@<1ex>[r]^{can} & {}^p\phi_t F\ar[l]^{var}
}
$$
Note that the objects in the diagram  are perverse sheaves on $\{0\}$
i.e.  vector spaces. Observe that since $T=I+var\circ can$ is the
monodromy of $F|_{\Delta^*}$, it is invertible.
This leads to the following elementary description of the category
 due to Deligne and Verdier (c.f. \cite[sect 4]{verdier}).

\begin{prop}\label{lemma:pervD}
  The category of perverse sheaves on the disk $\Delta$
which  are locally constant on $\Delta^*$ is
equivalent to the category of quivers of the form
$$
\xymatrix{
 \psi\ar@<1ex>[r]^{c} & \phi\ar[l]^{v}
}
$$
(i.e. finite dimensional vector spaces $\phi,\psi$ with maps as
indicated)
such that $I+v\circ c$ is invertible.
\end{prop}

It is instructive to consider some basic examples.
We see immediately that
$$
\xymatrix{
 0\ar@<1ex>[r] & V\ar[l]
}
$$
corresponds to the sky scraper sheaf $V_0$. 
Now suppose that $T$ is an automorphism of a vector space $V$.
This determines a  local system $L$ on $\Delta^*$.
Then the  perverse sheaf $j_*L[1]$ corresponds
to
$$
\xymatrix{
 V\ar@<1ex>[r]^{T-I} & im(T-I)\ar[l]^{v}
}
$$
where $v$ is the inclusion. The perverse sheaf $\R j_*L[1]$ 
corresponds to
$$
\xymatrix{
 V\ar@<1ex>[r]^{T-I} & V\ar[l]^{id}
}
$$

The above description can be extended to polydisks $\Delta^n$ \cite{galligo}. For
simplicity, we spell this out only for $n=2$. Let $t_i$ denote the coordinates.
Then we can attach to any perverse sheaf $F$, four vector
spaces $V_{11}={}^p\psi_{t_1}{}^p\psi_{t_2}F$,
$V_{12}={}^p\phi_{t_1}{}^p\psi_{t_2}F$... along with maps induced by $can$
and $var$. The set of these maps have additional commutivity and
invertibility constraints, 
such as $can_{t_1}\circ can_{t_2}= can_{t_2}\circ can_{t_1}$ etcetera.

\begin{thm}\label{thm:quivdelta2}
 The category of perverse sheaves on the polydisk $\Delta^2$
which  are constructible for the stratification $\Delta^2\supset \Delta\times
\{0\}\cup\{0\}\times \Delta\supset \{(0,0)\}$ is equivalent to the category
of  commutative quivers of the form
$$
\xymatrix{
 V_{11}\ar@<1ex>[r]^{c}\ar@<1ex>[d]^{c} & V_{12}\ar[l]^{v}\ar@<1ex>[d]^{c} \\ 
 V_{21}\ar[u]^{v}\ar@<1ex>[r]^{c} & V_{22}\ar[l]^{v}\ar[u]^{v}
}
$$ for which $I+v\circ c$ is invertible.
\end{thm}

It will be useful to characterize the subset of intersection cohomology complexes
among all the perverse ones. On $\Delta$ these are
either sky scraper sheaves $V_0$ in which case  $\phi=ker(v)$,
or sheaves of the form $j_*L[1]$ for which $\phi = im(c)$. On
$\Delta^2$ (and more generally), we have:

\begin{lemma}\label{lemma:varpluscan}
 A quiver
corresponds to  a direct sum of  intersection cohomology complexes
if and only if
$$V_{ij} = im(c)\oplus ker(v)$$
holds for every subdiagram
$$
\xymatrix{
 V_{i'j'}\ar@<1ex>[r]^{c} & V_{ij}\ar[l]^{v}
}
$$
\end{lemma}

\subsection{Kashiwara-Malgrange filtration}\label{section:KM}

 In this section, we give the $D$-module analogue of vanishing cycles.
But first we start with  a motivating example.

\begin{ex}\label{ex:KM}
  Let $Y=\C$  with coordinate $t$.  Let $M = \bigcup_m \OO_\C^n[t^{-m}]$
  with the standard basis $e_i$. We make this into a $D_\C$-module by letting
  $\dd_t$ act on the left through the singular connection
$\nabla=\frac{d}{dt}+A\frac{dt}{t}$, where  $A$ is a complex matrix with rational eigenvalues. 
There is no loss of generality in assuming that $A$ is in Jordan canonical form.
For each $\alpha\in \Q$,
define $V_\alpha M\subseteq M$ to be the $\C$-span of $\{t^{m}e_i\mid m\in \Z, m +a_{ii}\ge -\alpha\}$ where $e_i$ is the standard basis of $M$.
The following properties are easy to check:
\begin{enumerate}
\item The filtration is exhaustive  and left continuous:
$\cup V_\alpha M= M$ and
      $V_{\alpha+\epsilon} M = V_\alpha M$
for $0<\epsilon\ll 1$.
\item $ \partial_t V_\alpha  M\subseteq V_{\alpha+1}
      M$, and $tV_\alpha   M\subseteq V_{\alpha-1}   M$.
\item  The associated graded
$$Gr^V_\alpha  M = V_\alpha M/V_{\alpha-\epsilon} M
=  \bigoplus_{\{i\mid \alpha+a_{ii}\in \Z\}} \C t^{-\alpha-a_{ii}}e_i 
                $$
is isomorphic to the $-\alpha$ generalized eigenspace of $A$. Moreover,
the action of   $t\dd_t$ on this is space is identical to the action of $A$.
\end{enumerate}
(3) implies that the set of indices where  $V_\alpha  M$ jumps is
 discrete. (Such a filtration is called discrete.)
\end{ex}

Let $f:X\to \C$ be a holomorphic function, and
let $i:X\to X\times \C=Y$ be the inclusion
of the graph. Let $t$ be the coordinate on $\C$, and let
$$V_0 D_Y = \text{$D_{X\times\{ 0\}}$-module generated by }\{t^i\partial_t^j\mid i\ge j\}$$
This is the subring of differential operators preserving the ideal $(t)$.
Let $M$ be a regular holonomic $D_X$-module. It is called
 {\em quasiunipotent} along $X_0=f^{-1}(0)$ if   ${}^p\psi_f(DR(M))$ is quasiunipotent
with respect to the action of $T\in \pi_1(\C^*)$.
 Set $\tilde M = i_+M$. 

\begin{thm}[Kashiwara, Malgrange]
There exists at most one  filtration $V_\dt \tilde M$ on $\tilde M$ indexed by $\Q$,
such that
\begin{enumerate}
\item The filtration is exhaustive, discrete and left continuous.
\item  Each $V_\alpha \tilde M$ is a coherent $V_0D_Y$-submodule.
\item $ \partial_t V_\alpha \tilde M\subseteq V_{\alpha+1}\tilde
      M$,
and $tV_\alpha  \tilde M\subseteq V_{\alpha-1}  \tilde M$ with
equality for $\alpha< 0$.
\item $Gr^V_\alpha \tilde M$ is
a generalized eigenspace of $t\partial_t$ with eigenvalue $-\alpha$.
\end{enumerate}
 If $M$ is quasiunipotent along $X_0$, then $V_\dt\tilde M$ exists. 
\end{thm}

The module $\tilde M$ rather than $M$ is used to avoid problems caused by the singularities of $f$.
When $f$ is smooth, as in example~\ref{ex:KM}, we can and usually will construct the
$V$-filtration directly on $M$.

Given a perverse sheaf $L$ and $\lambda\in \C$, let ${}^p\psi_{f,\lambda} L$
and ${}^p\phi_{f,\lambda} L$ be the generalized $\lambda$-eigensheaves of ${}^p\psi_f L$ and ${}^p\phi_f
L$ under the $T$-action. The maps $can$ and $var$ respect this
decomposition into eigensheaves. It is technically convenient at this point
to define a modification $Var$ of $var$ such that $can\circ Var$ and
$Var\circ can$ equals $N=\frac{1}{2\pi \sqrt{-1}} \log T_u$, where $T_u$
is unipotent part of $T$ in the multiplicative Jordan
decomposition. After observing that  on each eigensheaf
$$\log T_u = \frac{1}{\lambda}N'
-\frac{1}{2}(\frac{1}{\lambda}N')^2+\ldots$$
where $N'= T-\lambda$. A formula for $Var$ is easily found by setting 
$$Var = a_0var+a_1var\circ N'+a_2var \circ (N')^2\ldots,$$
and solving for the coefficients.

\begin{ex}
Continuing with example \ref{ex:KM}. Let $L=DR(M)$. Then
${}^p\psi_t L= \ker e^*\nabla$, where $t=e(\tau)=\exp(2\pi \sqrt{-1}\tau)$.
A basis for this is given by the columns of the  fundamental solution of $e^*\nabla=0$ 
which is  $z= \exp(-2\pi \sqrt{-1}A \tau)$.  This 
has monodromy given by $T=\exp(-2\pi\sqrt{-1} A)$. Then $Gr^V_\alpha M$ is
the $-\alpha$ generalized eigenspace  of $A$, and this is isomorphic to the $\lambda$
generalized eigenspace of $T$ where $\lambda=e(\alpha)$.
\end{ex}

\begin{ex}
In the previous example, when $n=1,A=0$, the perverse sheaf corresponding to
$M$ is $L=\R j_* \C_{\C^*}$. We can see that $Gr^V_1  \tilde M=\C t^{-1}\cong {}^p\phi_t L$.
For the submodule $\OO_\C\subset M$, $DR(\OO)= \C$ and $Gr^V_1  \OO= 0 \cong {}^p\phi_t \C$.
\end{ex}

\begin{thm}[Kashiwara, Malgrange]\label{thm:kashmalg}
  Suppose that $L  = DR(M)$.  Let $\alpha\in \Q$ and $\lambda =
  e(\alpha)$. Then
  $$
DR(Gr^V_\alpha  \tilde M) =
\begin{cases}
  {}^p\psi_{f,\lambda} L & \text{ if $\alpha\in [0,1)$}\\
{}^p\phi_{f,\lambda} L &\text{ if $\alpha\in (0, 1] $}
\end{cases}
$$
The endomorphisms $t\partial_t+\alpha,\partial_t,t$ on the left corresponds to 
 $N, can, Var$ respectively, on the right.
\end{thm}

A proof is given in \cite[\S 3.4]{saito1}. Note that some care is needed  in comparing formulas here and in [loc. cit.]. There are shifts in indices according to whether one is working with left or right $D$-modules,
and decreasing or increasing  $V$-filtrations.

\section{ Pre-Saito Hodge theory}
\subsection{Hodge structures}
A (rational) {\em pure  Hodge structure} of weight $m\in \Z$ consists of a 
finite dimensional vector space $H_\Q$ with a bigrading
$$H= H_\Q\otimes \C = \bigoplus_{p+q=m} H^{pq}$$
satisfying $\bar H^{pq} = H^{qp}$. 
Such structures arise naturally from the cohomology of compact
K\"ahler manifolds. For smooth projective varieties, we have further
constraints namely the existence of a polarization on its cohomology.
A polarization on a weight $m$ Hodge structure $H$
is quadratic form $Q$ on $H_\Q$  satisfying the
Hodge-Riemann relations:
\begin{eqnarray}
  \label{eq:HR}
  Q(u,v) &=& (-1)^m Q(v,u)\\
Q(H^{pq},H^{p'q'}) &=& 0,\text{ unless } p=q',q=p'\\
\sqrt{-1}^{p-q} Q(u,\bar u) &>& 0,\text{ for } u\in H^{pq}, u\not=0
\end{eqnarray}

Given a Hodge structure of weight $m$, its Hodge filtration is
the decreasing filtration
$$F^pH = \bigoplus_{p'\ge p}H^{p',m-p'}$$
The decomposition can be recovered from the filtration by
$$H^{pq} = F^p\cap \bar F^{q}$$

Deligne extended Hodge theory to singular varieties. The key
definition is that of a  {\em mixed Hodge structure}. This consists 
of a   bifiltered vector space $(H, F, W)$, with $(H,W)$ defined over
$\Q$, such that $F$ induces a pure  Hodge
structure of weight $k$ on $Gr_k^WH$. Note that by tradition $F^\dt$ is
a decreasing filtration denoted with superscripts, while $W_\dt$ is
increasing. Whenever necessary, we can interchange increasing
and decreasing filtrations by $F_\dt= F^{-\dt}$. A pure Hodge
structure of weight $k$ can be regarded as a mixed Hodge structure
such that $Gr_k^W H = H$. A mixed Hodge structure is {\em polarizable}
if each $Gr_k^W H$ admits a polarization.
Mixed Hodge  structures  form a
category in the obvious way. 
Morphisms are rational linear maps preserving filtrations. 

\begin{thm}[Deligne]
  The category  of mixed Hodge structures is a $\Q$-linear abelian
  category with tensor product. The functors $Gr_k^W$ and $Gr_F^p$ are exact.
\end{thm}

The main examples are provided by:

\begin{thm}[Deligne]
The singular rational 
cohomology of a complex algebraic variety carries a canonical polarizable mixed Hodge structure.
\end{thm}

Finally, recall that there is a unique one dimensional pure Hodge structure $\Q(j)$ of
weight $-2j$ up to isomorphism. By convention, the lattice, is taken
to be $(2\pi \sqrt{-1})^j\Q$. The $j$th Tate twist $H(j)= H\otimes \Q(j)$.

\subsection{Variations of Hodge structures}
Griffiths introduced the notion a  {\em  variation of Hodge structure}  to describe
the cohomology of family of varieties $y\mapsto H^m(X_y)$, where
$f:X\to Y$ is a smooth projective map.
A variation of Hodge structure
of weight $m$ on a complex
manifold $Y$  consists of the following data:
\begin{itemize}
\item A locally constant sheaf $L$ of $\Q$ vector spaces with finite 
dimensional stalks.
\item A vector bundle with an integrable connection
$(E, \nabla)$ plus an isomorphism $DR(E)\cong  L\otimes \C[\dim Y]$.
\item A filtration $F^\dt$ of $E$ by subbundles satisfying 
Griffiths' transversality: $\nabla( F^p)\subseteq F^{p-1}$.
\item The data induces a pure  Hodge structure of weight $m$
on each of the stalks $L_y$.
\end{itemize}
A polarization on a variation of Hodge structure is a flat pairing $Q:L\times L\to \Q$ inducing
polarizations on the stalks. The main example is:

\begin{ex}\label{ex:gVHS}
  If $f:X\to Y$ is smooth and projective, $L = R^mf_*\Q$ underlies a
 polarizable variation of Hodge structure of weight $m$. $E=\R^m
 f_*\Omega_{X/Y}^\dt\cong \OO_Y\otimes L$ is the associated 
vector bundle with its Gauss-Manin connection and Hodge filtration
$F^p=\R^m f_*\Omega_{X/Y}^{\ge p}$.
\end{ex}

The next class of examples is of fundamental importance.

\begin{ex}\label{ex:nilorb}
An $n$ dimensional nilpotent orbit of weight $m$ consists of 
  \begin{itemize}
\item  a filtered complex vector space $(H,F)$ where $H$ carries a
  real structure,
\item a real form $Q$ on $H$ satisfying the first two Hodge-Riemman relations,
 \item $n$ mutually commutating nilpotent  infinitesimal isometries
   $N_i$ of $(H,Q )$ such that $N_iF^p\subseteq F^{p-1}$, 
\item   there exist $c>0$ such that for $Im(z_i)>c$, $(H,\exp(\sum
  z_iN_i)F,Q)$ is a polarized Hodge structure of
  weight $m$
  \end{itemize}
 Then the filtered polarized bundle  with fibre $(H,\exp(\sum
  z_iN_i)F,Q)$ descends to 
 a real holomorphic variation of Hodge structure on a  punctured polydisk
 $(\Delta-\{0\})^n$. The monodromy around the $i$th axis is precisely $\exp(N_i)$.
\end{ex}

Schmid \cite[\S 4]{schmid} showed that  any polarized variation of Hodge
structure on a punctured polydisk can be very well approximated by a nilpotent orbit
with the same asymptotic
behaviour. Thus the local study of a variation of Hodge structure can be reduced to nilpotent orbits.
When the local monodromies are unipotent (which can always be achieved by passing
to a branched covering), an approximating nilpotent orbit $(H, F,N_i)$ can  be constructed
in canonical fashion: $H$ is the fibre of the Deligne's canonical extension \cite{deligne-D} at $0$.
The Hodge bundles extend to subbundles of the canonical extension and $F$ is their fibre at $0$.
$N_i$ are the logarithms of local monodromy. 

The key analytic fact which makes the rest of the story possible is the following theorem
(originally proved by Zucker for curves).

\begin{thm}[Cattani-Kaplan-Schmid, Kashiwara-Kawai]\label{thm:IH}
Let $X$ be the complement of a divisor with normal crossings in a
compact K\"ahler manifold.
Then intersection cohomology with coefficients in a polarized variation of Hodge structure on $X$
is isomorphic to $L^2$ cohomology for a suitable complete K\"ahler metric on the
$X$.
\end{thm}

It should be noted that in contrast to the compact case,
 $L^2$ cohomology of a noncompact manifold is highly sensitive to the
 choice the metric, and that it could be infinite dimensional.
Here the metric is chosen to be asymptotic to a Poincar\'e metric along
transverse slices to the boundary divisor. From the theorem,
it follows that  $L^2$ cohomology is finite dimensional in this case.
When combined with the K\"ahler identities, we get

\begin{cor}
  Intersection cohomology $IH^i(X,L) = H^{i-n}(X,j_{!*}L[n])$ 
with coefficients in a polarized variation of Hodge structure of
weight $k$ carries a natural pure Hodge structure of weight $i+k$.
\end{cor}

\subsection{Monodromy filtration}

We start with a result in linear algebra which is crucial for later
developments.

\begin{prop}[Jacobson, Morosov]\label{prop:mono} Fix an integer $m$.
Let $N$ be a nilpotent endomorphism of a finite dimensional
vector space $E$ over a field of characteristic $0$.
 Then the there is a unique filtration
$$0\subseteq M_{m-l}\subset\ldots M_m\subseteq \ldots M_{m+l}=E$$
called the monodromy filtration on $E$ centered at $m$ (in symbols
$M=M(N)[m]$), characterized by following properties:
\begin{enumerate}
\item $N(M_k)\subseteq M_{k-2}$
\item $N^k$ induces an isomorphism $Gr_{m+k}^M(E)\cong Gr_{m-k}^M(E)$ 
\end{enumerate} 
\end{prop}

 In applications $N$ comes from monodromy, which explains the
 terminology. The construction of $M$ involves kernels, images
and other linear algebra operations, and so it applies to 
any artinian $\Q$ or $\C$-linear abelian category such as 
the categories of perverse sheaves or regular holonomic modules. 
The last part of the proposition is reminiscent of
the hard Lefschetz theorem.
There is an analogous decomposition into primitive parts: 

\begin{cor} 
  $$Gr^M_k(E) = \bigoplus_i N^i PGr^M_{k+2i}(E)$$
where  
$$PGr^M_{m+k}(E) = ker[N^{k+1}:Gr^M_{m+k}(E)\to Gr^M_{m-k-2}(E)]$$
\end{cor} 

Returning to Hodge theory, $M$ gives the weight filtration of
Schmid's limit mixed Hodge structure \cite{schmid},
which is a natural mixed Hodge structure on 
the space of nearby cycles. We combine this with the converse
statement of  Cattani and Kaplan \cite{ck}.

\begin{thm}[Schmid, Cattani-Kaplan]\label{thm:limitMHS}
Suppose that  $(H, F, N_i, Q)$ is a collection satisfying all  but the last
 condition of example \ref{ex:nilorb}. Then this is a nilpotent
orbit of weight $m$ if and only all of the following hold:
\begin{enumerate}
\item $M=M(\sum N_i)[m]=M(\sum t_iN_i)[m]$ for any $t_i>0$.
\item $(H, M)$ is a  real mixed Hodge structure.
\item This structure is polarizable; more specifically, the form
$Q(\cdot,N^k\cdot)$ gives a
polarization on the primitive Hodge structures $PGr^M_{m+k}(H)$.
\end{enumerate}

\end{thm}

Given a filtered vector space $(V,W)$ with nilpotent endomorphism $N$ preserving $W$,
a {\em relative monodromy } filtration is a filtration $M$ on $V$ such that
\begin{itemize}
\item $NM_k\subseteq M_{k-2}$
\item $M$ induces the monodromy filtration centered at $k$ on $Gr_k^WV$ with respect to $N$.
\end{itemize}

It is known that $M$ is unique if it exists (c.f. \cite{steen}),
although it need not exist
in general. In geometric situations, the existence of relative weight
filtrations was first established by Deligne using $\ell$-adic methods.

\subsection{Variations of mixed Hodge structures}

 A {\em variation of
mixed Hodge structure} consists of:

\begin{itemize}
\item A locally constant sheaf $L$ of $\Q$ vector spaces with finite 
dimensional stalks.
\item An ascending filtration $W\subseteq L$ by locally constant subsheaves
\item A vector bundle with an integrable connection
$(E, \nabla)$ plus an isomorphism $DR(E)\cong  L\otimes \C[\dim Y]$.
\item A filtration $F^\dt$ of $E$ by subbundles satisfying 
Griffiths' transversality: $\nabla( F^p)\subseteq F^{p-1}$.
\item $(Gr_m^W(L), \OO_X\otimes Gr_m^W(L), F^\dt(\OO_X\otimes
  Gr_m^W(L)))$ is a variation of pure Hodge structure of weight $m$.
\end{itemize}

 It follows that this data induces a mixed Hodge structure  on each of the stalks
 $L_y$. Steenbrink and Zucker \cite{steen} showed that additional
conditions are required to get a good theory. While these conditions
are rather technical, they do hold in most natural examples. 

A  variation of mixed Hodge structure over a
punctured disk $\Delta^*$ is {\em admissible} if

 \begin{itemize}
\item The pure variations $Gr_m^W(L)$ are polarizable.

\item There exists a limit Hodge filtration $\lim_{t\to 0} F_t^p$ compatible with
the one on $Gr_m^W(L)$ constructed by Schmid.

\item There exist a relative monodromy filtration $M$ 
on $(E=L_t,W)$ with respect to the logarithm  $N$ of the unipotent part of monodromy.
  \end{itemize}

For a general base $X$ the above conditions are required to  hold
for every restriction to a punctured disk \cite{kashiwara1}.
Note that pure polarizable variations of Hodge structure are automatically
admissible by \cite{schmid}. A local model for admissible
variations is provided by a {\em mixed nilpotent orbit} called an
infinitesimal mixed Hodge module in  \cite{kashiwara1}. This consists of:

\begin{itemize}
\item A bifiltered complex vector space $(H, F,W)$ with $(H,W)$
  defined over $\R$,
\item commuting real nilpotent endomorphisms $N_i$ satisfying
  $N_iF^p\subseteq F^{p-1}$ and $N_iW_k\subseteq W_k$, 
\end{itemize}
such that
\begin{itemize}
\item $Gr_k^W$ is a nilpotent orbit of weight $k$ for some choice of
  form.
\item The monodromy $M$ filtration exists for any partial sum $\sum
  N_{i_j}$, and $N_{i_j}M_k\subseteq M_{k-2}$ holds  for each ${i_j}$.
\end{itemize}

With an obvious notion of morphism, mixed nilpotent orbits of a given
dimension forms a category which turns out to be abelian \cite[5.2.6]{kashiwara1}.

\section{Hodge modules}
\subsection{Hodge modules on a curve}\label{subsection:MHcurve}

We are now ready to begin describing Saito category of Hodge modules
in a special case. We start with the  category $MF_{rh}(D_X, \Q)$
whose objects consist of
\begin{itemize}
\item A perverse sheaf $L$ over $\Q$.
\item A regular holonomic $D_X$-module $M$ with an isomorphism $DR(M) \cong L\otimes
      \C$.
\item A good filtration $F$ on $M$.
\end{itemize}
Morphisms are compatible pairs of morphisms of $(M,F)$ and $L$.
Variations of Hodge structure give examples of such objects (note that
Griffiths' transversality $\dd_{x_i}F_p\subseteq F_{p+1}$ is needed to
verify (3)).
The category $MF_{rh}(X,\Q)$ is really too big to do Hodge theory, and
Saito defines the subcategory of (polarizable)
Hodge modules  which provides a good setting. 
 The definition of this subcategory is extremely delicate,  so as a
 warm up we will give a direct  description for Hodge modules on  a 
smooth projective curve $X$ (fixed for the remainder of this section).

Given an inclusion of a point $i:\{x\}\to X$ and a polarizable pure Hodge structure
$(H,F, H_\Q)$ of weight $k$, the $D$-module pushforward $i_+H$ with the filtration induced
by $F$, together with the skyscraper sheaf $H_{\Q,x}$ defines an object of $MF_{rh}(X)$.
Let us call these polarizable Hodge modules of type $0$ and weight $k$.

Given a  polarizable variation of Hodge structure $(E,F,L)$ of weight $k-1$ over a Zariski 
open subset $j:U\to  X$, we define an object of $MF_{rh}(X)$ as follows.
The underlying  perverse sheaf is $j_*L[1]$ which is the intersection cohomology
complex for $L$. 
Since  $(E,\nabla)$ has regular singularities \cite{schmid}, 
we can extend it  to a vector bundle with logarithmic
connection on $X$ --  in many ways \cite{deligne-D}. The ambiguity depends on the
eigenvalues of the residues of the extension which 
are determined mod $\Z$. For every half open
interval $I$ of length $1$, there is a unique extension
$\bar E^{I}$ with eigenvalues in $I$. Let $\M' = \bigcup_I \bar E^I\subseteq j_*E$.
This is a $D_X$-module which corresponds to the perverse sheaf
$\R j_* L[1]\otimes \C$. Let $\M\subseteq \M'$ be the sub $D_X$-module generated
by $\bar E^{(-1,0]}$. This corresponds to what we want, namely $j_*L[1]$.
We filter this by 
$$F_p\bar E^{(-1,0]} = j_*F_pE\cap \bar E^{(-1,0]}$$
$$F_p \M =  \sum_i F_iD_XF_{p-i}\bar E^{(-1,0]}$$
Since the monodromy is  quasi-unipotent \cite[4.5]{schmid}, the $V$-filtration exists for general reasons.
However, in this case we can realize this explicitly by
$$V_\alpha \M = \bar E^{[-\alpha, -\alpha+1)}$$
This is essentially how we described it in example \ref{ex:KM}.
From this it follows that
$$F_p \cap V_\alpha \M = j_*F_pE\cap \bar E^{[-\alpha,-\alpha+1)}$$
for $\alpha < 1$.
The collection $(\M, F_\dt M, j_*L[1])$ defines an object of $MF_{rh}(X)$ that we
call a polarizable Hodge module of type $1$ of weight $k$. A polarizable Hodge module
of weight $k$ is a finite direct sum of objects of these two types. Let $MH(X,k)^{p}$ denote
the full subcategory of these.

\begin{thm}
$MH(X,k)^{p}$ is abelian and semisimple.
\end{thm}

In outline, the proof can be reduced to the following observations.
We claim that there are no nonzero morphisms between objects of type
$0$ and type $1$. To see this, we can replace $X$ by a disk and
assume that $x$ and $U$ above correspond to $0$ and $\Delta^*$ respectively.
Then by lemma \ref{lemma:pervD}, the perverse sheaves of type $0$ and $1$
 correspond to the quivers
$$
\xymatrix{
 0\ar@<1ex>[r] & \phi\ar[l]
}
$$
$$
\xymatrix{
 \psi\ar@<1ex>[r]^{onto} & \phi\ar[l]^{1-1}
}
$$
and the claim follows. We are thus reduced to dealing with the types
separately. For type $0$ (respectively $1$), 
we immediately reduce it the corresponding statements for the 
categories of  polarizable Hodge structures (respectively variation of
Hodge structures), where it is standard.
In essence the polarizations allow one to take orthogonal complements,
and hence conclude semisimplicity.

\begin{thm}
  If $\mathcal{M}\in MH(X,k)^{p}$, then its cohomology  $H^i(\mathcal{M})$
carries a pure Hodge structure of weight $k$.
\end{thm}

\begin{proof}
 Since $MH(X,k)^{p}$ is semisimple, we
can assume that $\mathcal{M}$ is simple. Then either $\mathcal{M}$ is supported
at point or it is of type $1$.
In the first case, $H^0(\mathcal{M}) = \mathcal{M}$ is already a Hodge structure 
of weight $k$ by definition, and the higher cohomologies vanish. 
In the second case, we appeal to theorem~\ref{thm:IH} or just
the special case due to Zucker.
\end{proof}

We end with a local analysis of Hodge modules. 
By a {\em Hodge quiver} we mean a diagram
$$
\xymatrix{
 \psi\ar@<1ex>[r]^{c} & \phi\ar[l]^{v}
}
$$
where $\psi$ and $\phi$ are mixed Hodge structures, $c$ is morphism
and $v$ is a morphism  $\phi\to \psi(-1)$ to the Tate twist, and such
that the compositions $c\circ v, v\circ c$ (both denoted by $N$) are
nilpotent. Actually, we will need to add a bit more structure, but we
can hold off on this for the moment.
We can encode the local structure of a Hodge module by a Hodge quiver.
A module of type $0$ and weight $k$ corresponds to the quiver
$$
\xymatrix{
 0\ar@<1ex>[r] & \phi\ar[l]
}
$$
where $\phi$ is pure of weight $k$. 

Next,  consider the type $1$ modules.
We wish to 
 associate a Hodge quiver to a polarizable variation of Hodge structure $(E,F,L)$ of weight $k-1$ 
over $\Delta^*$. First suppose that the monodromy is unipotent. Then let $\psi={}^p\psi_t L[1]$
with its limit mixed
Hodge structure i.e. the mixed Hodge structure on $\bar E^{(-1,0]}_0$ with Hodge filtration given by
$F^\dt \bar E^{(-1,0]}_0$ and weight
filtration $W= M(N)[k-1]$. Set $\phi =im(N)= \psi/ker(N)$ with its
induced mixed Hodge structure. Then 
$$
\xymatrix{
 \psi\ar@<1ex>[r]^{N} & \phi\ar[l]^{id}
}
$$
gives  the corresponding Hodge quiver. The Hodge quivers arising from this construction
are quite special in that $(\psi, F, N)$ and
$(\phi, F, N)$ are
nilpotent orbits of weight $k-1$ and $k$ respectively  by
\cite[2.1.5]{kk}. 

In general, the monodromy is quasi-unipotent. Thus
for some $r$, the pullback of $L$ along $\pi:t\mapsto t^r$ is
unipotent. So we can repeat the previous construction. But now the Galois
group $\Z/r\Z$ will act on everything, and it will be necessary to
include this action as part of the structure of a Hodge quiver so as not to lose any information.
 This can be understood in a directly without doing a base change.
The decomposition of $\psi={}^p\psi_t L[1]\cong {}^p\psi_{t}\pi^*L$
into a sum isotypic $\Z/r\Z$-modules can be identified with the
decompostion into generalized eigenspaces for the
$T$-action.  By
theorem \ref{thm:kashmalg}, we can realize this decomposition by
  $$\psi = \bigoplus_{0\le \alpha <1} Gr^V_\alpha \M.$$
For the Hodge filtration, we use the induced filtration
$$F^\dt Gr^V\M= \bigoplus_\alpha\frac{F^\dt \cap \ V_\alpha\M}
{F^\dt\cap   V_{\alpha-\epsilon}\M }
= \bigoplus_\alpha\frac{j_*F^\dt\cap \bar E^{[-\alpha,-\alpha+1)}}{j_*F^\dt\cap \bar E^{(-\alpha,-\alpha+1]}}$$
where $\M= D_\Delta \bar E^{(-1,0]}$ is defined as above. 
For $N$ we take the logarithm of the unipotent part of $T$ under the
multiplicative Jordan decomposition. Then $(\psi, F, N)$
is a nilpotent orbit of weight $k-1$, which determines a mixed Hodge
structure.  $\phi$ and the rest of the diagram is  defined as before.
In this case, we include the grading of $\psi,\phi$ by eigenspaces as part of the
structure of a Hodge quiver. In this way, we can recover the full
monodromy $T$ rather than just $N$.

In summary,  every polarizable Hodge module on the disk gives rise to a Hodge quiver.
This yields a functor which is necessarily faithful, since the quiver
determines the underlying perverse sheaf. This can be made into a full
embedding by  adding more structure to a Hodge quiver (namely,
 a choice of a variation of Hodge structure of
 $\Delta^*$...). However, this will not be necessary for our purposes.

\subsection{Hodge modules: overview} 

In the next section, we will define the  full subcategories 
$MH(X,n)\subseteq MF_{rh}(X,\Q)$ of {\em Hodge
modules} of weight $n\in \Z$ in general. Since this is rather
technical, we start by explaining the main results.

\begin{thm}[Saito]
  $MH(X,n)$ is abelian, and its objects possess strict support decompositions,
i.e.  that the maximal
sub/quotient module with  support in a given $Z\subseteq X$ can be split off as a direct
summand.
 There is an abelian subcategory $MH(X,n)^{p}$ of polarizable objects which is semisimple.
\end{thm}

 We essentially checked these properties for polarizable Hodge modules on curves in
section~\ref{subsection:MHcurve}. They have strict support decompositions
by the way we defined them.
Let $MH_Z(X,n)\subseteq MH(X,n)$ denote the subcategory of Hodge modules with 
strict support in $Z$, i.e. that all sub/quotient modules have support exactly $Z$.
The main examples are provided by the following.

\begin{thm}[Saito]\label{thm:MHvsVHS}
Any weight $m$ polarizable variation of Hodge structure $(L,\ldots)$
over an open subset of a closed subset 
$$U\stackrel{j}{\to} Z\stackrel{i}{\to} X$$
can be extended to a polarizable Hodge module in $MH_Z(X,n)^{p}\subseteq MH(X,n)^{p}$
with $n = m+ \dim Z$. 
The underlying perverse sheaf of the extension is the associated
intersection cohomology complex $i_*j_{!*}L[dim U]$. All simple
objects of   $MH(X,n)^{p}$ are of this form. 
\end{thm}

 Tate twists $\M\mapsto \M(j)$ can be
defined in this setting and are functors $MH(X,n)\to MH(X,n-2j)$.
Given a morphism $f:X\to Y$ and $L\in Perv(X)$, we can define perverse direct images
by ${}^pRf^i_* L = {}^p\mathcal{H}^i(\R f_* L)$. This operation extends
to Hodge modules. If $f$ is smooth and projective and $( M,F,L)\in
MH(X)$, the direct image is ${}^pRf^i_* L$ with  the associated
filtered $D$-module given by a Gauss-Manin construction:
$$ (\R^i f_*(\Omega_{X/Y}^\dt\otimes_{\OO_X} M ),\text{im} [\R^i f_*(\Omega_{X/Y}^\dt\otimes_{\OO_X}
F_{p+\dt}M )])$$
When applied to $(\OO_X, \Q_X[\dim X])$ with trivial filtration, we
recover the variation of Hodge structure of example \ref{ex:gVHS}
(after adjusting notation).

\begin{thm}[Saito]\label{thm:dirimageMH}
Let $f:X\to Y$ be a projective morphism with a relatively
ample line bundle $\ell $.
If  $\mathcal{M}=(M,F,L)\in MH(X,n)$ is polarizable, then
$${}^pR^if_*\mathcal{M}  \in MH(Y,n+i)$$
 Moreover, a hard Lefschetz theorem holds:
$$\ell^j\cup:{}^pR^{-j}f_*\mathcal{M}\cong {}^pR^{j}f_*\mathcal{M}(j)$$
\end{thm}

\begin{cor}
  Given a polarizable variation of Hodge structure defined
  on an open subset of $X$,  its intersection cohomology carries a
pure Hodge structure. This cohomology satisfies the
Hard  Lefschetz theorem.
\end{cor}

The last statement was originally obtained in the geometric
case in \cite{bbd}. The above results yield a Hodge theoretic
proof of the decomposition theorem of [loc. cit.]

\begin{cor}
  With assumptions of the theorem $\R f_*\Q$ decomposes into
a direct sum of shifts of intersection cohomology complexes.
\end{cor}

\subsection{Hodge modules: conclusion}

We now give the precise definition of Hodge modules. This is
given  by induction on dimension of the support.
This inductive process is handled via vanishing cycles. 
We start by explaining how to extend the construction to  $MF_{rh}(X,\Q)$.
Given a morphism $f:X\to \C$, and a $D_X$-module, $M$, 
we introduced the Kashiwara-Malgrange filtration $V$ on $M$
in section \ref{section:KM}. 
Now suppose that we have a good  filtration $F$ on $M$.
The pair $(M,F)$ is said to be be {\em quasi-unipotent and regular}
along $f^{-1}(0)$ if $V$ exists and if
\begin{enumerate} 
\item $t(F_pV_\alpha  \tilde M) = F_pV_{\alpha-1}  \tilde M$ for $\alpha < 1$.
\item $\partial_t(F_pGr^V_\alpha \tilde M) = F_{p+1} 
      Gr^V_{\alpha+1} \tilde  M$ for $\alpha\ge 0$.
\end{enumerate} 
hold along with certain finiteness properties that we will not spell out.
Saito \cite[3.4.12]{saito1} has shown theorem \ref{thm:kashmalg} holds assuming
only the existence of  such a filtration $F$. As for examples, note that
a variation of Hodge structure on a disk can be seen to satisfy these conditions  
using the formulas of  section  \ref{subsection:MHcurve} and \cite[3.2.2]{saito1}. 
Also any $D$-module which is quasi-unipotent and regular in the usual
sense admits  a filtration $F$ as above.

We extend  the functors $\phi$ and $\psi$ to $MF_{rh}(X,\Q)$ by
$$\psi_{f,e(\alpha)}(M,F, L) = (Gr^V_\alpha( \tilde M),F[1],
{}^p\psi_{f,e(\alpha)} L),\quad 0\le \alpha < 1$$
$$\phi_{f,e(\alpha)}(M,F, L) = ( Gr^V_{\alpha+1}( \tilde M),F,
{}^p\phi_{f,e(\alpha)} L),\quad 0\le \alpha < 1$$
$$\psi_f= \bigoplus \psi_{f,\lambda},\>\phi_f= \bigoplus \phi_{f,\lambda}$$
Where $F$ on the right  denotes the induced filtration. Thanks to the shift in $F$,
$\dd_t$ and $t$ induce morphisms
$$can: \psi_{f,e(\alpha)}(M,F, L) \to \phi_{f,e(\alpha)}(M,F, L) $$
$$var:\phi_{f,e(\alpha)}(M,F, L) \to \psi_{f,e(\alpha)}(M,F, L) $$
respectively.

We are ready to give the inductive definition of $MH(X,n)$.
The elementary definition for curves given earlier will turn out to be equivalent.
Define $X\mapsto MH(X,n)$ to be  the smallest collection of 
full subcategories of $MF_{rh}(X,\Q)$ satisfying:

\begin{enumerate} 
\item[(MH1)] If $(M,F, L)\in  MF_{rh}(X,\Q)$ has zero dimensional support, then it lies in
$MH(X, n)$ iff its stalks are Hodge structures of weight $n$. 

\item[(MH2)] If $(M,F,L)\in MHS(X,n)$ and $f:U\to \C$ is a general morphism from  
a Zariski open $U\subseteq X$, then  
\begin{enumerate}
\item $(M,F)|_U$ is quasi-unipotent and regular 
with respect to $f$.
\item $(M,F,L)|_U$ decomposes into a direct sum of a 
      module supported in $f^{-1}(0)$ and a module for which no
sub or quotient module is supported in $f^{-1}(0)$. 
    
\item If $W$ is the monodromy filtration of  $\psi_f(M,F,L)|_U$ (with respect to
the log of the unipotent part of monodromy) centered 
      at $n-1$, then
 $Gr^W_i\psi_f(M,F,L)|_U\in MH(U,i)$. Likewise for 
      $Gr^W_i\phi_{f,1}(M,F,L)$ with $W$ centered at $n$.
\end{enumerate} 
\end{enumerate}

This is a lot to absorb, so let us make few remarks about the definition.
\begin{itemize} 
\item If $f^{-1}(0)$ is in general position with respect to
$supp\, M$, the dimension of the support drops after applying the
functors $\phi$ and $\psi$. Thus this {\em is} an inductive definition.

\item  
The somewhat technical condition (b)  ensures that Hodge
modules admit strict support decompositions.
The condition can be rephrased as saying that $(M,F,L)$ splits as a sum
of the image of $can$ and the kernel of $var$. A refinement of
lemma~\ref{lemma:varpluscan} shows that $L$ will then decompose into a
direct sum of intersection cohomology complexes. This is ultimately
needed to be able to invoke theorem~\ref{thm:IH} when the time
comes to construct a Hodge structure on cohomology.

\item Although $MF_{rh}(X,\Q)$ is not an abelian category, the
category of compatible pairs consisting of a $D$-module and perverse
sheaf is. Thus we do get  a $W$ filtration for $\psi_f(M,F,L)|_U$ in (c)
by  proposition~\ref{prop:mono} by first suppressing $F$,
and then using the induced filtration.

\end{itemize}

There is a notion of polarization in this setting.
Given $(M,F,L)\in MH_Z(X,n)$, a polarization is a pairing
$S:L\otimes L\to \Q_X[2\dim X](-n)$ satisfying certain axioms.
The key conditions are  again inductive.  
When $Z$ is a point, $S$ should correspond to a polarization 
on the Hodge structure at the stalk in the usual sense.
 In general, given a (germ of a) function $f:Z\to \C$ which is not identically
zero, $S$ should induce a polarization on the nearby cycles $Gr_\dt^W \psi_fL$
(using the same recipe as in theorem \ref{thm:limitMHS} (3)).
Once all the  definitions are in place, the proofs of the theorems involve a
rather elaborate induction on dimension of supports.

In order to get a better sense of what Hodge modules look like, we
extend the local description given in \S\ref{subsection:MHcurve}  to 
the polydisk $\Delta^2$.  We assume that the
underlying perverse sheaves satisfying the assumptions of theorem \ref{thm:quivdelta2}.
 Then using a vanishing cycle construction refining the one used earlier,
such a polarizable Hodge module gives rise to a two dimensional
Hodge quiver. This  is  a commutative diagram of $\Q^2$-graded mixed Hodge structures
$$
\xymatrix{
 H_{11}\ar@<1ex>[r]^{c}\ar@<1ex>[d]^{c} & H_{12}\ar[l]^{v}\ar@<1ex>[d]^{c} \\ 
 H_{21}\ar[u]^{v}\ar@<1ex>[r]^{c} & H_{22}\ar[l]^{v}\ar[u]^{v}
}
$$ 
such that the maps labeled $c$ and $v$ are required to
satisfy the same conditions as  in \S\ref{subsection:MHcurve}.
In particular, compositions $c\circ v$, $v\circ v$  are nilpotent.
These nilpotent transformations are denoted by $N_1$ for  horizontal arrows,
and $N_2$ for vertical.

Given a Hodge module $M$ of weight $k$,
 $H_{11} = \psi_{t_2}\psi_{t_2}M\in MF_{rh}(pt,\Q)$ etcetera are
 filtered vector spaces with a $\Q^2$-grading by eigenspaces.
 These filtrations  are used as the Hodge filtrations.
The weight filtration is given by the shifted monodromy filtration with respect to $N=N_1+N_2$,
where $N_i$ are the logarithms of the unipotent parts of monodromy around the axes. The shifts
are indicated in the diagram below
\begin{equation}\label{eq:weights}
\xymatrix{
 k-2\ar@<1ex>[r]\ar@<1ex>[d] & k-1\ar[l]\ar@<1ex>[d] \\ 
 k-1\ar[u]\ar@<1ex>[r] & k\ar[l]\ar[u]
}
\end{equation}
 These quivers are quite special. A Hodge quiver is
  pure of weight
$k$ if the conditions of lemma~\ref{lemma:varpluscan} hold, and if
each $(H_{ij}, F, W, N_1,N_2)$ arises from a  pure  nilpotent orbit of
weight  indicated indicated in \eqref{eq:weights}.
With this description in hand, we can verify a local version of theorem \ref{thm:MHvsVHS}.
Given a polarizable variation of Hodge structure $L$ of weight $k-2$ on $(\Delta-\{0\})^2$ 
corresponding to a nilpotent orbit $(H,F, N_1,N_2)$, we can construct the
Hodge quiver
$$
\xymatrix{
 H\ar@<1ex>[r]^{N_1}\ar@<1ex>[d]^{N_2} & N_1H\ar[l]^{id}\ar@<1ex>[d]^{N_2} \\ 
 N_2H\ar[u]^{id}\ar@<1ex>[r]^{N_1} & N_1N_2H\ar[l]^{id}\ar[u]^{id}
}
$$ 
where $H$ is equipped with the monodromy filtration $W(N_1+N_2)$ on the upper left.
The remaining vertices are given the mixed Hodge structures induced on the images of this.
The fact that these vertices correspond to nilpotent orbits of the correct weight follows
from   \cite[1.16] {cks} or \cite[2.1.5]{kk}. Note that $j_{!*}L[2]$
is  the perverse sheaf corresponding to this diagram.
The details can be found in  \cite[\S 3]{saito2}.

\subsection{Mixed Hodge modules}

Saito has  given an extension of the previous theory
by defining the notion of mixed Hodge module.
We will define a {\em pre-mixed Hodge module} on $X$ to consist of

\begin{itemize}
\item A perverse sheaf $L$ defined over $\Q$, together with a
 filtration $W$ of $L$ by perverse subsheaves.
\item A regular holonomic $D_X$-module $\M$ with a filtration $W\M$ 
which corresponds to $(L\otimes \C, W\otimes \C)$ under
Riemann-Hilbert.
\item A good filtration $F$ on $\M$.
\end{itemize}

These objects form a category, and Saito defines the
subcategory of {\em mixed Hodge modules} $MHM(X)$ by a rather delicate induction.
The key points are that for $(\M,F,L,W)$ to be in $MHM(X)$, he requires
\begin{itemize}
\item The associated graded objects $Gr_k^W(\M,F,L)$
yield polarizable Hodge modules of weight $k$. 
\item  For any (germ of a) function $f$ on $X$,
the relative monodromy filtration $M$ (resp. $M'$) for
$\psi_f(\M,F,L)$ (resp. $\phi_{f,1}(\M,F,L)$) with respect to $W$ exists.
\item The pre-mixed Hodge modules
 $(\psi_f(\M,F,L),M)$ and  $(\phi_{f,1}(\M,F,L),M')$ are in fact mixed Hodge modules
on $f^{-1}(0)$.
\end{itemize}
This is not a complete list of all the conditions. Saito also requires
that mixed Hodge modules should be extendible across divisors in a
sense we will not attempt to make precise.

 The main properties are summarized below:
 
\begin{thm}[Saito]
\begin{enumerate}
\item[]
\item $MHM(X)$ is abelian, and it contains each $MH(X,n)^{p}$ as a full
       abelian subcategory.
\item $MHM(\text{point}) $ is the category of polarizable mixed
 Hodge structures.
\item If $U\subseteq X$ is open, then any admissible variation of mixed 
 Hodge structure on $U$ extends to an object in $MHM(X)$.
 \end{enumerate}
\end{thm}

 Finally, we have:

\begin{thm}[Saito]
  There is a realization functor of triangulated categories
 $$real:D^bMHM(X)\to D^b_{constr}(X,\Q)$$
and functors
$$ {}^p\cH^i: D^b MHM(X)\to MHM(X)$$
$$ \R f_*: D^b MHM(X)\to D^b MHM(Y)$$
$$ \mathbb{L} f^*: D^bMHM(Y)\to D^bMHM(X)$$
for each morphism $f:X\to Y$, such that
$${}^p\cH^i\circ real= real\circ {}^p\cH^i$$
$$real( \R f_* \mathcal{M}) = \R f_*real(\mathcal{M})$$
$$real( \mathbb{L} f^* \mathcal{N}) = \mathbb{L} f^*real(\mathcal{M})$$
Similar statements hold for various other standard operations such
as tensor products and cohomology with proper support.
\end{thm}

Putting these results together yields  Saito's mixed Hodge structure:

\begin{thm}\label{thm:saitoMHS}
  The cohomology of a smooth variety $U$ with coefficients in an
admissible variation of mixed Hodge structure $L$ carries a canonical
mixed Hodge structure. 
\end{thm}

When the base is a curve, this was first proved by   Steenbrink and
Zucker \cite{steen}.

\begin{proof}
$L$ gives an object of $MHM(U)$. Thus $H^i(\R p_*L)$ carries a mixed Hodge structure,
where  $p:U\to pt$ is the projection to a point. 
\end{proof}

We extend the local description of perverse sheaves and Hodge modules to 
mixed Hodge modules.
 A one dimensional {\em filtered Hodge quiver} is a diagram
$$
\xymatrix{
 \psi\ar@<1ex>[r]^{c} & \phi\ar[l]^{v}
}
$$
where $\psi$ and $\phi$ are filtered $\Q$-graded mixed Hodge structures, $c$ and $v:\phi\to \psi(-1)$
are morphisms  and  the compositions $c\circ v, v\circ c$ (both denoted by $N$) are nilpotent.
So each object $\psi,\phi$ is equipped with $3$ filtrations $F,M,W$, where $M$ denotes the
weight filtration. For a filtered Hodge quiver to arise  from a mixed Hodge module,
we require that 
\begin{itemize}
\item $Gr_W^k$ is pure of weight $k$.
\item $(\psi, F, W, N)$ and $(\phi, F, W, N)$ are mixed nilpotent orbits with
$M$ the relative monodromy filtration. 
\end{itemize}
We call these {\em admissibility} conditions.
This extends to polydisks in a straightforward manner.

\subsection{Explicit construction}\label{section:expconstr}

We will  give a bit more detail on the construction of the mixed Hodge
structure in theorem~\ref{thm:saitoMHS}. 
For a different perspective see \cite{elzein}.
Let $U$ be a smooth $n$ dimensional variety.
We can choose a smooth compactification $j:U\to X$ such that $D=X-U$ is
a divisor with normal crossings. Fix an admissible variation of
mixed Hodge structure $(L,W,E, F,\nabla)$ on $U$. 
Extend $(E,\nabla)$ to a vector bundle $\bar E^I$ with a logarithmic
connection such that the eigenvalues of its residues lie in $I$
as in section \ref{subsection:MHcurve} Let $\bar E=\bar E^{(-1,0]}$
and $\M = \bigcup_I \bar E^I\subseteq j_*E$.
Then $\M$  is a $D_X$-module which corresponds to the perverse sheaf
$\R j_* L[n]\otimes \C$. Filter this by 
$$F_p\M =\sum_i F_iD_XF_{p-i}\bar E$$
Then:

\begin{thm}
  There exists compatible filtrations $\tilde W$ on $\R j_* L[n]$ and $\M$ extending $W$
over $U$  such that 
$(\R j_* L[n],  \tilde W, \M, F)$ becomes a mixed Hodge module.
\end{thm}

We describe $\tilde W$ when  $n=\dim X\le 2$.
We start with the case of $n=1$, which is
essentially due to Steenbrink and Zucker \cite{steen}. 
 It suffices to construct $\tilde W$ locally near any point of $p\in D$. Choose a coordinate disk
$\Delta$ at $p$. The perverse sheaf $\R j_* L[1]$
corresponds to the quiver
\begin{equation}
  \label{eq:quiv1}
  \xymatrix{
  \psi\ar@<1ex>[r]^{N} & \psi\ar[l]^{id}
}
\end{equation}
where $\psi=\psi_t L$ and 
$N$ is the logarithm of unipotent part of monodromy at $p$. 
 We  make this into an admissible filtered
Hodge quiver by replacing $\psi$ on the left with the mixed
nilpotent orbit $(H, F^\dt H, W_\dt H, N)$, where we can take $H=\bar E_0$
when the monodromy is unipotent, and $H=\oplus Gr^V_\alpha (D_\Delta \bar
E)$ in general. On the
right, we keep $(H, F^\dt H, N)$ but replace $W$ with 
\begin{equation}
  \label{eq:Nstar}
(N_*W)_{k-1} := N W_{k} + M_k\cap W_{k-1}  
\end{equation}
where $M$ is the relative monodromy filtration of $(H,W)$. This is a
mixed nilpotent orbit by \cite[5.5.4]{kashiwara1}.
Then define the filtration $\tilde W_k$ of (\ref{eq:quiv1}) by
$$\xymatrix{
  W_{k-1}\ar@<1ex>[r]^{N} & (N_*W)_{k-2}\ar[l]^{id}
}
$$
The key point is  that $Gr_k^{\tilde W}$ is pure of weight $k$,
because  it is a sum  of
$$
\xymatrix{
  Gr^W_{k-1}\ar@<1ex>[r]^{N} & N Gr^W_{k-1}\ar[l]^{id}
}
$$
and 
$$
\xymatrix{
  0\ar@<1ex>[r] & Gr^M_{k-2}(W_{k-2}/NW_{k-2})(-1)\ar[l]
}
$$
both pure of weight $k$. A proof can be extracted from
 \cite[pp 514-516]{steen} or \cite{kashiwara1}. Since we have
an equivalence of categories between perverse sheaves on $\Delta$ and quivers,
$\tilde W$ above determines a filtration of $\R j_* L[1]$ by perverse subsheaves.

We now turn to the surface case. Replace $X$ by a polydisk $\Delta^2$
such that $D$ corresponds to (a subdivisor of) the union of axes. An
admissible filtered Hodge quiver is now a commutative diagram
$$
\xymatrix{
 H_{11}\ar@<1ex>[r]^{c}\ar@<1ex>[d]^{c} & H_{12}\ar[l]^{v}\ar@<1ex>[d]^{c} \\ 
 H_{21}\ar[u]^{v}\ar@<1ex>[r]^{c} & H_{22}\ar[l]^{v}\ar[u]^{v}
}
$$ 
where $H_{ij}$ arise from mixed nilpotent orbits. The maps labeled $c$ and $v$ are required to
satisfy the same conditions as above. 
Let $N_i$ denote the logarithms of the unipotent parts of local monodromies of $L$ around
the axes. Then the perverse  sheaf $\R j_* L[2]$ is represented by the
quiver
$$
\xymatrix{
 H\ar@<1ex>[r]^{N_1}\ar@<1ex>[d]^{N_2} & H\ar[l]^{id}\ar@<1ex>[d]^{N_2} \\ 
 H\ar[u]^{id}\ar@<1ex>[r]^{N_1} & H\ar[l]^{id}\ar[u]^{id}
}
$$ 
where $H=\bar E_{0}$ in the unipotent case.
 We make this into a filtered Hodge quiver by equipping this
with $F^\dt H,N_i$ and the weight filtrations given below
$$
\xymatrix{
 W_k\ar@<1ex>[r]\ar@<1ex>[d] & (N_{1*}W)_{k-1}\ar[l]\ar@<1ex>[d] \\ 
 (N_{2*}W)_{k-1}\ar[u]\ar@<1ex>[r] & (N_{1*}N_{2*}W)_{k-2}\ar[l]\ar[u]
}
$$
where $N_{1*}W$... are defined as in  \eqref{eq:Nstar}. Note that
symmetry holds for the last filtration $N_{1*}N_{2*}W = N_{2*}N_{1*}W$
by \cite[5.5.5]{kashiwara1}.
Then $\tilde W_k= W_{k-2}$ gives the desired filtration on $\R j_* L[2]$.
 Again arguments of \cite{kashiwara1} can  be
used to verify that $Gr^{\tilde W}_k$ is pure of weight $k$, as required.
\bigskip

Finally, let us describe Saito's  mixed Hodge structure explicitly
by elaborating on a remark in \cite{saito1.5}. We have an isomorphism
$$\alpha:\R j_* L\otimes \C\cong \Omega_X^\dt(\log D)\otimes \bar E^{[0,1)}$$
and the filtrations defining the Hodge structure can be displayed
rather explicitly using the complex on the right.
We define two filtrations $F$ and $\tilde W$ on this complex. 
$$F^p(\Omega_X^\dt(\log D)\otimes \bar E^{[0,1)}) = \bigoplus_i
\Omega_X^i(\log D)\otimes F^{p-i}\bar E^{[0,1)}$$
We construct $\tilde W$ by taking the filtration induced by 
$(\Omega_X^\dt\otimes \tilde W_{\dt+n} \M)[n]$
under the inclusion
$$\Omega_X^\dt(\log D)\otimes \bar E^{[0,1)}\subset
\Omega_X^\dt\otimes \M[n]$$
We recall at this point that Deligne  introduced a device,
 called  a {\em cohomological mixed Hodge complex} \cite[\S 8.1]{deligne-H}, for
producing mixed Hodge structures. This consists of 
bifiltered complex $(A_\C,W_\C, F)$ of sheaves of  $\C$-vector spaces,
a filtered complex of $(A,W)$ of sheaves over $\Q$, and
a filtered quasi-isomorphism $(A,W)\otimes \C\cong (A_\C,W_\C)$.
The crucial  axioms are that
\begin{itemize}
\item  this datum should induce a pure weight $i+k$ Hodge
  structure on $H^i(Gr^W_kA)$
\item the filtration induced by $F$ on $Gr^W_k$ is strict, i.e.
the map $H^i(F^pGr^W_kA)\to H^i(Gr^W_kA)$ is injective.
\end{itemize}

\begin{prop} With these filtrations
  $$(\R j_*L, \tilde W; \Omega_X(\log D)^\dt\otimes  \bar E, \tilde W, F;\alpha)$$
 becomes a cohomological mixed Hodge complex.
\end{prop}

To verify the above axioms, one observes that $Gr^{\tilde W}_k\R j_*L$ 
decomposes into a direct
sum of intersection cohomology complexes associated to pure variations
of Hodge structure  of the correct weight. An appeal to
theorem~\ref{thm:IH} shows that $H^i(Gr_k^{\tilde W}\R j_*L)$ carries
pure Hodge structures of weight  $k+i$. Strictness  follows from
similar considerations.

From \cite[\S 8.1]{deligne-H}, we obtain

\begin{cor}
The Hodge filtration is induced by $F$ under the isomorphism
$$H^i(U,L\otimes \C)\cong \mathbb{H}^i(\Omega_X^\dt(\log D)\otimes
\bar E^{[0,1)})$$
The weight filtration is given by
$$W_{i+k} H^i(U,L\otimes \C) = \mathbb{H}^i(X, W_k(\Omega_X^\dt(\log D)\otimes \bar E^{[0,1)})).$$ 
\end{cor}

Finally, we remark that in the simplest case $L=\Q$, $\tilde W$ above
coincides with the filtration $W$ defined by Deligne \cite[\S 3.1]{deligne-H}; in
particular,
Saito's Hodge structure coincides with Deligne's in this case.

\section{Comparison}
From now  on let $f:X\to Y$ be a smooth projective morphism of smooth quasiprojective
varieties. We recall two results

\begin{thm}[Deligne {\cite{deligne-L}}]\label{thm:delL}
The Leray spectral sequence
$$E_2^{pq} = H^p(Y,R^qf_*\Q) \Rightarrow H^{p+q}(X,Q)$$
degenerates. In particular $E_2^{pq} \cong Gr^p_LH^{p+q}(X)$
for the associated ``Leray filtration'' $L$.
\end{thm}

\begin{thm}[{\cite{arapura}}] There exists  varieties $Y_p$ and morphisms
$Y_p\to Y$ such that
$$L^pH^i(X,\Q) = ker[H^i(X,\Q)\to H^i(X_p,\Q)]$$
where $X_p = f^{-1}Y_p$.
\end{thm}

It follows that each $L^p$ is a filtration by
sub mixed  Hodge structures. When combined theorem \ref{thm:delL},
we get a mixed Hodge structure
on $H^p(Y,R^qf_*\Q) $ which we will call the
naive mixed Hodge structure. On the other hand,
$R^qf_*\Q$ carries a pure hence admissible variation
of Hodge structure, so we can apply Saito's theorem \ref{thm:saitoMHS}.

\begin{thm}
  The naive mixed Hodge structure on $H^p(Y,R^qf_*\Q) $  coincides with Saito's.
\end{thm}

\begin{proof}
 Note that $R^if_*\Q$ are
local systems and hence perverse sheaves up to shift.
More specifically,
$$R^{i}f_*\Q = real({}^p\cH^{i-\dim X+\dim Y}\R f_*\Q[\dim X])[-\dim Y]$$
  Deligne \cite{deligne-L} actually proved a stronger version of the above
theorem which implies that 
\begin{equation}\label{eq:lefdecomp}
\R f_*\Q \cong \bigoplus_i R^if_*\Q[-i]
\end{equation}
(non canonically) in $D_{constr}^b(X)$.
By theorem~\ref{thm:dirimageMH} and [loc. cit.], we have the corresponding decomposition
$$\R f_*\Q[\dim X] = \bigoplus_i {}^p\cH^i(\R f_*\Q[\dim X])[-i]$$
 in  $D^bMHM(Y)$. Note that the Leray filtration
is induced by the truncation filtration 
$$
L^pH^i(X,\Q) = image[H^i(Y,\tau_{\le i-p} \R f_*\Q\to H^i(Y,\R 
f_*\Q)]
$$
Under (\ref{eq:lefdecomp}), 
$$\tau_{\le p} \R f_*\Q \cong \bigoplus_{i\le p} R^if_*\Q[-i]$$
Therefore we have an
isomorphism of mixed Hodge structures
$$H^{i}(X,\Q) \cong \bigoplus_{p+q=i} H^p(Y,R^pf_*\Q)$$
 where the right side is equipped with Saito's mixed
Hodge structure. Under this isomorphism, $L^p$ maps to
$$H^{i-p}(R^pf_*\Q)\oplus H^{i-p+1}(R^{p-1}f_*\Q)\oplus\ldots$$
The theorem now follows.
\end{proof}

This theorem would follow also from the remark \cite[4.6.2]{saito2}
that one can build a $t$-structure on $D^bMHM(X)$ which agrees with the
standard one on $D^b_{constr}(X)$. However, we decided to present an alternate proof,
since a detailed justification for the remark has not been given.

\begin{appendix}

\section{Supplement to \cite{arapura}}

We will clarify the proof of \cite[lemma 3.13]{arapura}, which
  contained a hidden assumption about the behaviour of
filtered acyclic resolutions under d\'ecalage.
I will freely use the  notation from \cite{arapura}.
To simplify terminology, assume that
complexes are always bounded below, and filtrations are biregular.
Given a sheaf $\F$ on a space $X$, let $\G^\dt(\F)$ denote its  canonical
flasque resolution, constructed by Godement \cite{godement}:
$\G^0(\F) = \prod_x \F_x$ etcetera. 
If $\F^\dt$ is  a  complex,
let $\G^\dt(\F^\dt)$ denote the total complex resulting from this construction.
We record the basic properties, which are either standard or easily checked.

\begin{lemma}
$\G^\dt$ gives an exact functor from complexes to complexes of flasque sheaves.
The canonical map  $\F^\dt\to \G^\dt(\F^\dt)$ is a quasi-isomorphism. If
 $\E$ is a sheaf on a closed set $i:Z\to X$, 
$i_*\G^\dt(\E) = \G^\dt(i_*\E)$. If  $\E,\F$ are sheaves on an open set $j:U\to X$
and on $X$ respectively, then 
$j_!\G^\dt(\E) = \G^\dt(j_!\E)$ and $j^*\G^\dt(\F) = \G^\dt(j^*\F)$. 
\end{lemma}

The following was already observed in \cite{deligne-H}.
 
\begin{cor}
  If $(\F^\dt, F^\dt)$ is a filtered complex, $((\G^\dt(\F^\dt), \G^\dt(F^\dt))$
is a filtered acyclic resolution. (Here acyclicity is with respect
direct image, and $\Gamma$ in particular.)
\end{cor}

\begin{cor}
  If $(\F^\dt, F^\dt)$ is a filtered complex, then
  $$((\G^\dt(\F^\dt), \G^\dt(Dec(F^\dt))) =((\G^\dt(\F^\dt), Dec(\G^\dt(F^\dt))) $$
is a filtered acyclic resolution of $(\F^\dt, Dec(F^\dt))$.
\end{cor}

\begin{proof}
  The fact that $Dec$ commutes with $\G^\dt$ follows from exactness of $\G^\dt$.
The rest follows from the previous corollary.
\end{proof}

The following should be viewed as a replacement for \cite[lemma 3.3]{arapura}.

\begin{cor}
  If $(Y, Y_\dt)$ is an object of $FV_\C$ and $\F$ is a sheaf on $Y$,
$$(\G^\dt(\F), S^\dt(Y_\dt, \G^\dt(\F)) = (\G^\dt(\F), \G^\dt(S^\dt(Y_\dt, \F))$$
is a filtered acyclic resolution of $(\F, (S^\dt(Y_\dt, \F))$.
\end{cor}

\begin{proof}
 The fact that $S^\dt$ commutes with $\G^\dt$ follows from the last part
of the lemma. The rest follows from corollary 1. 
\end{proof}

Note that a filtered quasi-isomorphism $(\A, F\A)\to (\B, F\B)$
induces a filtered quasi-isomorphism  $(\A, Dec(F\A))\to (\B, Dec(F\B))$.
(This can be deduced from \cite[Prop 1.3.4]{deligne-H}.)
Thus $Dec$ is well defined on the filtered derived category.

We will indicate the corrected proof of \cite[lemma 3.13]{arapura} below; referring
to the original notation and giving
only the modifications.

\begin{proof}[Proof of lemma 3.13]
  Let $(\I^\dt, \Sigma^\dt) = (\G^\dt(\F), S^\dt(X_\dt, \G^\dt(\F)))$.
This gives a filtered acyclic resolution of $\F$ with respect to both 
$ S(X_\dt,\F)$ and $Dec(S(X_\dt,\F))$ by the above discussion. 
As in \cite{arapura}, we obtain a map of filtered complexes
$$(f_*\I,\tau)\to (f_*\I, Dec(f_*\Sigma)) = (f_*\I, f_*(Dec(\Sigma)))$$
These  complexes are filtered acyclic, and the map  induces a morphism 
$$ \R f_*(\F, \tau)\to \R f_*(\F, Dec(S(X_\dt,\F)))$$
in the filtered derived category of sheaves on $Y$. Applying $\R \Gamma$
yields a morphism in the filtered derived category of abelian groups, which
results in a map of spectral sequences
$$
\begin{array}{ccc}
  E_1(\R\Gamma(\F,\tau))&\to &  E_1(\R\Gamma(\F, Dec(S(X_\dt, \F)))) \\ 
\| & & \| \\
 \LL_2(f,\F)& \to &  E_2(X_\dt, \F)
\end{array}
$$
as required. The vertical identifications and
the  proof of naturallity are the same as in \cite{arapura}.
\end{proof}

\end{appendix}


\end{document}